\documentclass[10pt]{amsart}
\usepackage{fancyhdr}
\usepackage{amssymb,amscd}
\usepackage{amsthm,amsmath}
\usepackage[cp850]{inputenc}

\renewcommand{\O}{\mathcal{O}_C}
\newcommand{\OJ}{\mathcal{O}_{J(C)}}
\renewcommand{\L}{\mathcal{L}}

\newcommand{\lra}{\longrightarrow}

\newcommand{\oc}{\omega_C}

\newcommand{\Z}{\mathbb Z}

\def\ord{\mathrm{ord}\,}
\def\mult{\mathrm{mult}\,}

\theoremstyle{plain}
\newtheorem{thm}{Theorem}[section]
\newtheorem{lem}[thm]{Lemma}
\newtheorem{prop}[thm]{Proposition}
\theoremstyle{remark}
\newtheorem{rem}{Remark}

\begin{document}

\thispagestyle{empty}

\title[The Vanishing of the Theta Function]
      {The Vanishing of the Theta Function in the KP Direction: a Geometric
         Approach}
\author[Ch. Birkenhake]{Christina Birkenhake}
  \thanks{Supported by DFG-contracts Hu 337/5-1}
  \address{Christina Birkenhake\\Mathematisches Institut\\
           Bismarkstra\ss e 1 1/2\\ D-91054 Erlangen}
  \email{birken@mi.uni-erlangen.de}
  \urladdr{http://www.mi.uni-erlangen.de/\~{}birken/}

\author[P. Vanhaecke]{Pol Vanhaecke}
  \address{Pol Vanhaecke\\Universit\'e de Poitiers, D\'epartement de
Poitiers,
           T\'el\'eport 2, Boulevard Marie et Pierre Curie, BP 30179,
F-86962
           Futuroscope Chasseneuil Cedex}
  \email{Pol.Vanhaecke@mathlabo.univ-poitiers.fr}
  \urladdr{http://wwwmathlabo.sp2mi.univ-poitiers.fr/\~{}vanhaeck/}

\renewcommand{\subjclassname}{\textup{2000} Mathematics Subject
     Classification}

\begin{abstract}
  We give a geometric proof of a formula, due to Segal and Wilson, which
  describes the order of vanishing of the Riemann theta function in the
  direction which corresponds to the direction of the tangent space of a
  Riemann surface at a marked point. While this formula appears in the work
  of Segal and Wilson as a by-product of some non-trivial constructions
  from the theory of integrable systems (loop groups, infinite-dimensional
  Grassmannians, tau functions, Schur polynomials, $\dots$) our proof only
  uses the classical theory of linear systems on Riemann surfaces.
\end{abstract}

\maketitle
\tableofcontents

\section{Introduction}\label{SW_intro}

The fundamental paper \cite{segalwilson} by Segal and Wilson on soliton
equations leads to an explicit formula for computing the vanishing of the
Riemann theta function in a direction which is natural from the geometric
point of view. In order to present this formula, let $C$ be a compact
Riemann surface (of genus $g>1$), let $p\in C$ and let $\Theta$ denote the
theta divisor $\Theta\subset\mathrm{Pic}^{g-1}(C)$. Also let $\vartheta$
denote Riemann's theta function, for which $(\vartheta)=\Theta$. For a
point $\L\in\Theta$, consider the embedding
$C\to\mathrm{Pic}^{g-1}(C):q\mapsto \L(q-p)$. The natural direction alluded
to above is the tangent space $X_p$ to this embedded curve at
$\L$. Following \cite{segalwilson} the vanishing of $\vartheta$ at $\L$ in
the direction of $X_p$, denoted $\ord_\L(\vartheta,X_p)$ is obtained by
considering the infinite subset of $\Z$, defined by
$$
  S_\L=\{s\in\mathbb Z\mid h^0(\L((s+1)p))=h^0(\L(sp))+1\}.
$$
In fact, the sections of $\L$ over $C\setminus\{p\}$ define an
infinite-dimensional plane $W$, which is an element of the Sato
Grassmannian and the vanishing of the tau function in the KP-direction is,
according to \cite[Prop.\ 8.6]{segalwilson}, given by the codimension of
$W$, which is explicitly given by the finite sum $\sum_{i\geq 0}i-s_i$, upon writing
$S_{\L}=\{s_0<s_1<s_2<\cdots\}$. The tau function coincides, up to an
exponential factor, with the Riemann theta function of $C$
(\cite[Th. 9.11]{segalwilson}) and the tangent direction $X_p$ coincides
with the KP-direction (see \cite[Lemma 5 and Appendix
0]{shiota}). Therefore, the order of vanishing is given by
$$
  \mathrm{ord}_\L(\vartheta,X_p)=\sum_{i\geq 0}i-s_i.
  \label{SWintro}
$$

The purpose of this paper is to give an algebraic-geometric proof of this
result. Our proof uses (only) the classical theory of linear systems on
Riemann surfaces and it highlights the geometric meaning of the order of
vanishing. As is pointed out in \cite[footnote p.\ 51]{segalwilson} an
independent (analytical) proof of this formula has also been given by John
Fay (see \cite{fay}), by using the theory of theta functions. 

The first step of our approach consists of an interpretation of the order of
vanishing as the intersection multiplicity of the theta divisor with a copy
of $C$, properly embedded (at least around $\L$) in
$\mathrm{Pic}^{g-1}(C)$. If we pull back the theta divisor using this
embedding we find a divisor $R$ on $C$ which is the sum of the ramification
divisors of the maps
\begin{equation*}
  \varphi_\L^k:C\lra\mathrm{Grass}_{k+1}(H^0(\L)^*),
\end{equation*}
which are the natural generalizations of the morphisms $\varphi_\L:C\lra
\mathbb{P}(H^0(\L)^*)$ defined by the linear system $\vert \L\vert$
(assumed here base point free). It follows that the order of vanishing is
given by the multiplicity of $p$ in $R$, leading to 
\begin{equation*}
  \mathrm{ord}_\L(\vartheta,X_p)=\sum_{i=1}^{n}m_i-i.
\end{equation*}
where $\{m_1<\cdots<m_n\}$ is the gap sequence $G_p(\L(np))$ of $\L(np)$ at
$p$. This formula is independent of $n$, which is assumed sufficiently
large (e.g. $n=g$ will do). Noticing that for $n=g$ one has $s_i=g-m_{g-i}$
(for $i=0,\dots,g-1$) from which Formula \eqref{SWintro} follows at once. 

Notice that the Segal-Wilson formula for the vanishing of the tau function
may also be applied in the case of tau functions that come from singular
curves. It would be interesting to adapt our geometric arguments to this
case, leading to a formula for the vanishing of the \emph{theta functions
for singular curves}, as proposed in \cite[Remark 6.13]{segalwilson}.

The structure of this paper is as follows. In Section 2 we fix the notation
and we recall the notions of gap numbers for arbitrary line bundles. In
Section 3 we translate the order of vanishing of the theta function in
terms of intersection theory and we show that this order is given as an
inflectionary weight. This is used in Section 5 to obtain an explicit
formula, which we show to be equivalent to the formula by Segal and
Wilson. 

\section{Preliminaries}\label{prel}
\par In this section we introduce the notation and collect some results on
curve theory. Throughout the whole paper $C$ denotes a compact Riemann
surface of genus $g$ and $p\in C$ a marked point.

For a divisor $D$ on $C$ we denote by $\O(D)$ the corresponding line bundle
and for a line bundle $\L$ on $C$ its linear system is denoted by $\vert
\L\vert$. We use the standard abbreviations $h^0(\L)$ for $\dim H^0(C,\L)$
and $\L(D)$ for $\L\otimes\O(D)$, where $\L$ is any line bundle and $D$ is
any divisor on $C$. We will use the Riemann-Roch theorem in the form
\begin{equation*}
  h^0(\L)=h^0(\oc\otimes \L^{-1})-g+\deg(\L)+1,
\end{equation*}
where $\L$ is any line bundle on $C$ and $\oc$ is the canonical bundle of
$C$.

\smallskip

We now recall the notions of gap numbers and inflectionary
weights. For proofs and details we refer to \cite[Sect.\ VII.4]{miranda}
and to \cite[Ch.\ 1 Ex.\ C]{acgh}.

Let $\L$ be a line bundle on $C$ of positive degree and let $q\in C$. An
integer $m\geq1$ is called a \emph{gap number for $\L$ at $q$} if
\begin{equation*}
  h^0(\L(-mq)) = h^0(\L(-(m-1)q))-1,
\end{equation*}
and the set $G_q(\L)$ of gap numbers for $\L$ at $q$ is called the
\emph{gap sequence of $\L$ at $q$}; its cardinality is $r=h^0(\L)$ and no
gap number is larger than $\deg \L+1$. Writing
\begin{equation*}
  G_q(\L)=\{ 1\leq m_1< m_2 < \cdots < m_{r}\leq \deg \L+1\},
\end{equation*}
we have that $m_1>1$ if and only if $q$ is a base point of $\L$ and that
$m_{r}=\deg \L+1$ if and only if $\L=\O(\deg \L\cdot q)$. For a general
point $q\in C$ the gap sequence of $\L$ at $q$ is $\{1,2,\ldots,h^0(\L)\}$;
a point $q$ for which the gap sequence of $\L$ at $q$ is not of this form
is called an \emph{inflection point for $\L$}.  Notice that $q$ is an
inflection point if and only if $h^0(\L(-rq))\not=0$, where $r=h^0(\L)$.

If the linear system $\vert \L\vert$ is base point free the inflection
points have the following geometric interpretation. Consider the morphism
$\varphi_\L:C\lra \mathbb{P}(H^0(\L)^*)$ defined by the linear system
$\vert \L\vert$. For a generic $q\in C$ there is a unique $k$-dimensional
osculating plane to $\varphi_\L(C)$ at $\varphi_\L(q)$, yielding a
well-defined morphism
\begin{equation}
  \varphi_\L^k:C\lra\mathrm{Grass}_{k+1}(H^0(\L)^*),
\end{equation}
called the $k$-th \emph{associated map}. This way one arrives at
$h^0(\L)-1$ associated maps $\varphi_\L^{i-1},\ i=1,\dots,h^0(\L)-1,$
($\varphi_\L^0=\varphi_\L$). In these terms a point $q$ is an inflection
point if and only if $q$ is a ramification point of one of the maps
$\varphi_\L^k$. We denote the ramification divisor of $\varphi_\L^k$ by
$R_k(\L)$ and we define
\begin{equation*}
        R(\L)=\sum_{k=1}^{h^0(\L)-1} R_{k-1}(\L).
\end{equation*}
The multiplicity $w_q(\L)$ of $q$ in $R(\L)$ is called the
\emph{inflectionary weight} of $q$ with respect to $\L$ and is given by
\begin{equation}\label{inflec_weight}
   w_q(\L) = \sum_{i=1}^{h^0(\L)}(m_i-i).
\end{equation}
%
%
When $\L$ is not base point free we define the inflectionary weights
$w_q(\L)$ by \eqref{inflec_weight} and the ramification divisor $R(\L)$ by
$R(\L)=\sum_qw_q(\L)q$. This divisor admits an alternative description as
the zero divisor of $W=W(z)(dz)^{\frac{n(n-1)}{2}}$ where $z$ is a local
coordinate, $n=h^0(\L)$ and $W(z)=W(f_1,\ldots,f_n)$ is the Wronskian with
respect to any basis $f_1,\ldots,f_n$ of $H^0(\L)$. In particular $W$ is a
holomorphic section of the line bundle $\L^n\otimes\oc^{\frac{n(n-1)}{2}}$.

Taking $\L=\oc$ one recovers the well-known notion of the gap sequence
of $q\in C$, denoted by $G_q$, and the above definition of the
inflectionary points and weights reduces, by a simple application of
Riemann-Roch, to the standard definition of Weierstra\ss\ points and
their weights.

\smallskip

Finally let us fix our conventions about the Jacobian $J(C)$ of
$C$. By definiton $J(C)=H^0(\oc)^*/H_1(C,\Z)$, so that the vector space 
$H^0(\oc)^*$ is canonically identified with the tangent space of $J(C)$ 
at every point. By the Abel-Jacobi theorem there is a canonical isomorphism 
$J(C)\simeq\mathrm{Pic}^{0}(C)$. Moreover every line bundle $\L$ on $C$ of degree $g-1$ 
induces an isomorphism $J(C)\simeq\mathrm{Pic}^{0}(C) \lra \mathrm{Pic}^{g-1}(C),
\,\mathcal P\mapsto \L\otimes \mathcal P$. For our purposes it is convenient 
to work with $\mathrm{Pic}^{g-1}(C)$ rather than $\mathrm{Pic}^{0}(C)$. 
So in the sequel we identify $J(C)$ with $\mathrm{Pic}^{g-1}(C)$ 
without further noticing, the underlying isomorphism 
(respectively line bundle defining the isomorphism) will always be evident 
from the context. The main advantage working with 
$J(C)=\mathrm{Pic}^{g-1}(C)$ is, that we have a canonical  
\emph{theta divisor} $\Theta=\{\L\in J(C)\mid h^0(\L)>0\}$.
By Riemann-Roch, $\Theta$ is invariant with
respect to the natural involution
\begin{equation}
  \iota : J(C) \lra J(C), \, \iota(\L) = \oc\otimes \L^{-1}.
\end{equation}
More precisely we have $h^0(\L)=h^0(\iota(\L))$ for any $\L\in\Theta$.

We denote by $\vartheta$ the
Riemann theta function on $H^0(\oc)^*$ for which $\pi^*\Theta$ is the zero
divisor of $\vartheta$, where $\pi$ is the natural projection
$H^0(\oc)^*\to J(C)$. 

For any $\L\in J(C)$ we have an embedding
$\alpha_{\L,p}$ of $C$ into $J(C)$, given by
$\alpha_{\L,p}(q)=\L(q-p)$. Clearly for different $\L$ and $p$ the maps
$\alpha_{\L,p}$ only differ by a translation on $J(C)$.

\section{Geometric description of the order of vanishing}
\par Let $\L\in J(C)$ and let $X$ be a one-dimensional
subvector space of the tangent space $H^0(\oc)^*$ at $\L$. Choose any point $l$ in
the fiber of $\pi$ over $\L$ and consider the affine line $l+X$ which
passes through $l$ and which has direction~$X$. The order of vanishing of
$\vartheta\vert_{l+X}$ at the point $l$ is independent of the choice of
$l\in\pi^{-1}(\L)$. So define the \textit{order of vanishing of}
$\vartheta$ \textit{at} $\L$ \textit{in the direction of} $X$, denoted
$\ord_\L(\vartheta,X)$, as $\ord_l\vartheta\vert_{l+X}$.  If $\Theta$ does
not contain the straight line $\bar X=\pi(l+X)$ then there exists a small
neighborhood $U$ of $l$ in $l+X$ such that $\pi(U)\cap\Theta=\{\L\}$ and
$\ord_\L(\vartheta,X)=(\pi(U)\cdot\Theta)_\L$, the intersection
multiplicity of $\Theta$ with $\pi(U)$ at $\L$.

Let $X_p$ denote the tangent space to $\alpha_{\L,p}(C)$ at $\L$. Notice
that, as a subvector space of $H^0(\oc)^*$, $X_p$ does not depend on $\L$ but only
on the point $p\in C$. We wish to compute $\ord_\L(\vartheta,X_p)$ for an
arbitrary $\L\in\Theta$; if $\L\notin\Theta$ then this order is trivially
zero. It can be shown \footnote{This follows from \cite[Prop. 8.6 and Th.\
9.11]{segalwilson} but a geometric proof of this (geometric!) property is
unknown.} that for all $C$, $\L$ and $p$ that $\bar X_p$ is not contained in
$\Theta$, which is clearly true as soon as $C$, $\L$ or $p$ is generic. For
this the idea is to replace $\pi(U)$ by a a complete curve which, around
$\L$, looks like $\pi(U)$.  Notice that if $\L\in\Theta$ we cannot use
$\alpha_{\L,p}(C)$ because the latter does not necessarily intersect
$\Theta$ properly. Consider for any integer $n\neq 0$ the morphism
\begin{equation}
   \alpha_{\L,p,n} :C\lra J(C),\,\, \alpha_{\L,p,n}(q)=\L(nq-np).
\end{equation}
Notice that $\alpha_{\L,p,n}(p)=\L$ and that for a small neighborhood
$V$ of $p$ in $C$ the tangent space to $\alpha_{\L,p,n}(V)$ at $\L$ is
precisely $X_p$.
\begin{lem}\label{lem31} For all $\L\in J(C)$ and $n>0$ we have
\begin{enumerate}
\item  $\alpha_{\L,p,n}(C)$ and $\Theta$ intersect properly if and only if
    $h^0(\L(-np))=0$;
\item $\alpha_{\L,p,-n}(C)$ and $\Theta$ intersect properly if and only if
     $h^0(\iota(\L)(-np))=0$.
  \end{enumerate}
\end{lem}
\begin{proof}
We prove (1), the proof of (2) is similar. Recall that an
irreducible curve intersects a divisor properly precisely when the curve is
not contained in the support of the divisor. So $\alpha_{\L,p,n}(C)$ and
$\Theta$ do not intersect properly if and only if $h^0(\L(nq-np))>0$ for all
$q\in C$. We claim that this is equivalent to $h^0(\L(-np))>0$. Indeed, by
Riemann-Roch
\begin{align*}
  h^0(\L(nq-np))&=h^0(\iota(\L)(np-nq)),\quad\text{and}\\
  h^0(\iota(\L)(np))&=h^0(\L(-np))+n.
\end{align*}
So $h^0(\L(np-nq))>0$ for all $q\in C$ if and only if
$h^0(\iota(\L)(np))>n$, leading to our claim.
\end{proof}
In particular, when $\vert n\vert\geq g$ then
\begin{equation*}
  \ord_\L(\vartheta,X_p)=(\Theta\cdot \alpha_{\L,p,n}(V))_\L,
\end{equation*}
where $V$ is a small neighborhood of $p$ in $C$. Pulling this intersection
back to $C$ we get that for any $\vert n\vert\geq g$
\begin{equation*}
  \ord_\L(\vartheta,X_p)=\mult_p(\alpha_{\L,p,n}^*\Theta).
\end{equation*}
This multiplicity will be computed in the next section.


\section{The divisor $\alpha_{\L,p,-n}^*\Theta$}

The aim of this section is to prove the following

\begin{thm}\label{thm42}
For all $\L \in J(C)$, $p\in C$, and $n>0$ with $h^0(\iota(\L)(-np))=0$
$$
\alpha_{\L,p,-n}^*\Theta=R(\L(np)).
$$
\end{thm}
For the proof we need the following
\begin{prop} \label{prop411}
For all $\L\in J(C), p\in C$ and $n>0$
$$
\alpha_{\L,p,-n}^*\OJ(\Theta)=(\L(np))^n\otimes\oc^{\frac{n(n-1)}{2}}.
$$
\end{prop}
\begin{proof}\
\par
\textsc{Step I:} The case $n=1$ follows exactly from \cite{cav} Lemma
11.3.4 with $x=0, \kappa=\L$, and $c=p$. \\

\textsc{Step II:} For $n\geq 1$ consider the difference map
\begin{equation*}
  \delta_\L^n:C^{2n}\longrightarrow J(C),\quad
  \delta_\L^n(p_1,q_1,\ldots,p_n,q_n)=\L(\textstyle\sum_ip_i-q_i)
\end{equation*}
and denote by $\pi_i^n:C^{2n}\longrightarrow C$ the $i$-th projection.
We show by induction on $n$ that for all $n\geq 1$ and all $\L\in J(C)$
\smallskip
\begin{equation}\label{eqn411}
\begin{split}
{\delta_\L^n}^*\OJ
(\Theta)&=\\
\bigotimes_{i=1}^n\Bigl({\pi_{2i-1}^n}^*(\oc\otimes\L^{-1})\otimes{\pi_{2i}^n}^*\L\Bigr)
&\otimes\mathcal O_{C^{2n}}\Bigl(\sum_{1\leq i<j\leq
2n}(-1)^{i+j+1}(\pi_i^n,\pi_j^n)^*\Delta\Bigr),
\end{split}
\end{equation}
where $\Delta$ denotes the diagonal in $C^{2}$.

For $n=1$ we have to show that 
\begin{equation*}
  {\delta_\L^1}^*\OJ(\Theta)
  ={\pi_1^1}^*(\oc\otimes\L^{-1})\otimes{\pi_2^1}^*\L\otimes\mathcal O_{C^2}(\Delta)
\end{equation*}
for all $\L\in \mathrm{Pic}^{g-1}(C)$. According to the
Seesaw Principle (see \cite{cav} A.9) it suffices to show that the
restrictions to $C\times\{q\}$ and $\{q\}\times C$ of both sides of the
equation coincide for all $q\in C$.  But since the composition of
$\delta_\L^1$ with the natural embedding $C\simeq
C\times\{q\}\longrightarrow C\times C$ is the map $\iota\circ
\alpha_{\oc\otimes\L^{-1},q,-1}$ and $\iota^*\Theta=\Theta$ we have, using
Step I,
\begin{align*}
  {\delta_\L^1}^*\OJ(\Theta)\vert C\times\{q\}
     &=\alpha_{\oc\otimes\L^{-1},q,-1}^*\OJ(\Theta)=\oc\otimes\L^{-1}(q)\\
     &={\pi_1^1}^*(\oc\otimes\L^{-1})\otimes{\pi_2^1}^*\L\otimes\mathcal O_{C^2}(\Delta)\vert C\times\{q\},
\end{align*}
and similarly for the restriction to $\{q\}\times C$.

Now suppose $n>1$ and equation \eqref{eqn411} holds for all $n'<n$. Restricting
both sides of equation \eqref{eqn411} to $C^{2n-2}\times\{p,q\}$ and
$\{p_1,q_1,\ldots,p_{n-1},q_{n-1}\}\times C^2$ for all
$p,q,p_1,q_1,\ldots,p_{n-1},q_{n-1}\in C$, and using the
induction hypothesis for $n'=n-1$ and $n'=1$ respectively, the Seesaw
Principle implies that also equation \eqref{eqn411} holds.
\\
\goodbreak
\textsc{Step III:} Consider the embedding $\jmath_p:C\longrightarrow
C^{2n},\quad \jmath_p(q)=(p,q,\ldots,p,q)$ and notice that ${\delta_\L^n}\circ
\jmath_p=\alpha_{\L,p,-n}$, so that
\begin{equation*}
  \alpha_{\L,p,-n}^*\OJ(\Theta)=
  \jmath_p^*{\delta_\L^n}^*\OJ(\Theta).
\end{equation*} 
It follows that $\alpha_{\L,p,-n}^*\OJ(\Theta)$ can be computed from
\eqref{eqn411}. Since
\begin{equation*}
  (\pi_{i}^n,\pi_{j}^n)\circ\jmath_p(q)=\left\{
  \begin{array}{cl}
    (p,p)&i,j\hbox{ odd}\\
    (q,q)&i,j\hbox{ even}\\
    (p,q)\hbox{ or } (q,p)&\hbox{otherwise,}
  \end{array}
  \right.
\end{equation*}
we have that 
\begin{equation*}
  \jmath_p^*(\pi_{i}^n,\pi_{j}^n)^*\mathcal O_{C^2}(\Delta)=\left\{
  \begin{array}{cl}
    \O&i,j\hbox{ odd}\\
    \oc^{-1}&i,j\hbox{ even}\\
    \O(p)&\hbox{otherwise.}
  \end{array}
  \right.
\end{equation*}
It follows that the pull back by $\jmath_p$ of the right hand side of
\eqref{eqn411} equals $\L^n(n^2p)\otimes\oc^{\frac{n(n-1)}{2}}$. This
completes the proof.
\end{proof}
\begin{proof}[Proof of Theorem \ref{thm42}]
By the choice of $n$ we have $h^0(\L(np))=n$ and the curve
$\alpha_{\L,p,-n}(C)$ intersects the divisor $\Theta$ properly. The line
bundle $(\L(np))^n\otimes\oc^{\frac{n(n-1)}{2}}$ has two distinguished
divisors, namely $\alpha_{\L,p,-n}^*\Theta$ (according to Proposition
\ref{prop411}) and $R(\L(np))$ (according to Section \ref{prel}). Moreover
these divisors have the same support, since by definition $q\in
\alpha_{\L,p,-n}^*\Theta$ if and only if $\alpha_{\L,p,-n}(q)\in \Theta$,
i.e., $h^0(\L(np-nq))>0$, and this is the case if and only if $q $ is an
inflection point for the line bundle $\L(np)$.
%
%
For generic $\L$ and $p$ the line bundle $\L(np)$ admits only normal
inflection points, so $R(\L(np))=\sum_{i=1}^{n^2g} q_i$ with pairwise
different points $q_i$, and hence $R(\L(np))=\alpha_{\L,p,-n}^*\Theta$.
This equality extends to all $\L$ and $p$ for which
$\alpha_{\L,p,-n}^*\Theta$ exists and by Lemma \ref{lem31} this is exactly
the set $\{(\L,p)\in J(C)\times C\,\vert\,h^0(\iota(\L)(-np))=0\}$.
\end{proof}

\begin{rem}
  Similarly one can show that if $\alpha_{\L,p,n}(C)$ intersects $\Theta$
  properly, then $\alpha_{\L,p,n}^*\Theta=R(\iota(\L)(np))$.
\end{rem}

\section{Formula(s) for the order of vanishing}
\par
In this section we prove the following theorem:

\begin{thm}\label{thm41}
For every $\L\in \Theta$
 \begin{equation}\label{eqn41}
\mathrm{ord}_\L(\vartheta,X_p)=\sum_m m-g+h^0(\L((g-m)p)),
\end{equation}
where the sum runs over the $g$ integers $m$ satisfying
$h^0(\L((g-m)p))=h^0(\L((g-m+1)p))-1$.
\end{thm}

We will obtain it as a direct consequence of the following proposition.

\begin{prop}\label{prop43}
With $n$ chosen such that $\alpha_{\L,p,-n}(C)$ intersects $\Theta$ properly
(e.g., $n=g$), the order of vanishing of $\vartheta$ at $\L$ in the
direction $X_p$ is the inflectionary weight of $p$ with
respect to $\L(np)$. Therefore,
\begin{equation}\label{eqn42}
  \mathrm{ord}_\L(\vartheta,X_p)=\sum_{i=1}^{n}m_i-i.
\end{equation}
where $\{m_1<\cdots<m_n\}$ is the gap sequence $G_p(\L(np))$ of $\L(np)$ at
$p$.
\end{prop}
\begin{proof}
According to Section 3 and Theorem \ref{thm42} we have that
\begin{equation*}
\ord_\L(\vartheta,X_p)=\mult_p(\alpha_{\L,p,-n}^*\Theta)=w_p(\L(np))=
\sum_{i=1}^{n}m_i-i.
\end{equation*}
\end{proof}

\begin{proof}[Proof of Theorem \ref{thm41}]


By definition the sum in equation \eqref{eqn41}
runs over the set $G_p(\L(gp))=\{m_1<\cdots<m_g\}$ of gap numbers of $\L(gp)$
at $p$. An immediate computation shows that  $h^0(\L((g-m_i)p))=g-i$ for
$i=1,\ldots,g$. So the assertion follows from Proposition \ref{prop43} with
$n=g$.
\end{proof}

We now relate Theorem \ref{thm41} to the Formula \eqref{SW_intro}, given by
Segal and Wilson. Recall from the introduction the infinite set
$$
  S_\L=\{s\in\mathbb Z\mid h^0(\L((s+1)p))=h^0(\L(sp))+1\}.
$$
\begin{prop}\label{prop45}
Denote $S_{\L}=\{s_0<s_1<s_2<\cdots\}$. Then
$$
\mathrm{ord}_\L(\vartheta,X_p)=\sum_{i\geq 0}i-s_i.
$$
\end{prop}
\begin{proof}
Note first that $s_0\geq -\deg\L-1=-g$. For $-g\leq s\leq g-1$ we have that
$s\in S_\L$ if and only if $g-s\in G_p(\L(gp))=\{m_1<\cdots<m_g\}$, so
that $s_i=g-m_{g-i}$ for $i=0,\dots,g-1$. On the other hand $n\in
S_{\L}$ for any $n\geq g$ so that $s_n=n$ for any $n\geq g$. Summing up
we find
\begin{equation*}
  \sum_{i\geq0}{i-s_i}=\sum_{i=0}^{g-1}i-g+m_{g-i}=\sum_{i=1}^gm_i-i.
\end{equation*}
Hence Formula \eqref{SW_intro} follows from Proposition \ref{prop43}.
\end{proof}

\end{document}